\documentclass[12pt]{article}
\usepackage{color}
\definecolor{darkblue}{rgb}{0.00,0.25,0.50}
\usepackage[colorlinks,filecolor=blue,citecolor=darkblue]{hyperref}

\usepackage{amsfonts,amssymb,amsmath,amsthm}
\usepackage{url}
\usepackage{enumerate}
\usepackage[ukrainian, russian, english]{babel}
\usepackage[cp1251]{inputenc}
\begin{document} \selectlanguage{ukrainian}
\thispagestyle{empty}

\title{}

УДК 517.51 \vskip 3mm

\begin{center}
\textbf{\Large Order estimates of approximative characteristics
\\ of functions from classes \boldmath{$S^{r}_{1,\theta}B(\mathbb{R}^d)$}}
\end{center}
\begin{center}
\textbf{\Large Порядкові оцінки апроксимативних характеристик \\  функцій з класів \boldmath{$S^{r}_{1,\theta}B(\mathbb{R}^d)$}}
\end{center}
\vskip0.5cm

\begin{center}
 S.~Ya.~Yanchenko\\ \emph{\small
Institute of Mathematics NAS of
Ukraine, Kyiv}
\end{center}
\begin{center}
C.~Я.~Янченко \\
\emph{\small Інститут математики НАН України, Київ}
\end{center}
\vskip0.5cm

\begin{abstract}

We obtained exact order estimates of approximation of the classes
$S^{\boldsymbol{r}}_{1,\theta}B$ by entire functions of exponential
type with supports of their Fourier transforms in step hyperbolic cross. The error of the approximation
estimated in the metric of Lebesgue spaces, $L_q(\mathbb{R}^d)$, $1<q\leqslant \infty$.
\vskip 3 mm

 Одержано точні за порядком оцінки наближення функцій  з класів
$S^{\boldsymbol{r}}_{1,\theta}B(\mathbb{R}^d)$ за допомогою цілих функцій експоненціального типу з носіями їх перетворення Фур'є у східчастому гіперболічному хресті. Похибка наближення оцінюється в метриці простору Лебега, $L_q(\mathbb{R}^d)$, $1<q\leqslant \infty$ .
\end{abstract}

\vskip 0.5 cm


У роботі продовжено вивчення апроксимативних характеристик
  класів функцій Нікольського--Бєсова
  $S^{\boldsymbol{r}}_{p,\theta}B(\mathbb{R}^d)$~\cite{Nikolsky_63}, \cite{Amanov_1965}
 у просторі $L_q(\mathbb{R}^d)$. Встановлено точні за
порядком оцінки наближення функцій із згаданих класів цілими функціями, зі
спектром зосередженим на множині, яка називається східчастим
гіперболічним хрестом. Основна увага приділяється випадку, коли $p=1$.

\vskip 3 mm  \textbf{1. Означення класів функцій та апроксимативних характеристик}.
Нехай $\mathbb{R}^d$~--- $d$-вимірний евклідів простір з елементами
${\boldsymbol{x}=(x_1,...,x_d)}$ і
${(\boldsymbol{x},\boldsymbol{y})=x_1y_1+...+x_dy_d}$.  Нехай
${L_q(\mathbb{R}^d)}$, ${1\leqslant q\leqslant\infty}$,
--- простір вимірних на $\mathbb{R}^d$ функцій ${f(\boldsymbol{x})=f(x_1,...,x_d)}$
зі скінченною нормою
 $$
\|f\|_q:=
\Bigg(\int\limits_{\mathbb{R}^{d}}|f(\boldsymbol{x})|^{q}d\boldsymbol{x}
\Bigg) ^{\frac{1}{q}}, \ 1\leqslant q<\infty,
 $$
 $$
  \|f\|_{\infty}:=\mathop {\rm ess \sup}\limits_{\boldsymbol{x}\in \mathbb{R}^d}
  |f(\boldsymbol{x})|.
 $$

Для функції $f(\boldsymbol{x})\in L_q(\mathbb{R}^d)$ визначимо
різницю $1$-го порядку з кроком $h$ за змінною $x_j$ таким чином:
$$
 \Delta_{h,j}f(\boldsymbol{x})=f(x_1,\dots,x_{j-1},x_j+h,x_{j+1},\dots,x_d)-f(\boldsymbol{x})
$$
 і, відповідно, $l$-го порядку, $l \in \mathbb{N}$,
$$
 \Delta_{h,j}^{l}f(\boldsymbol{x})=
 \overbrace{\Delta_{h,j}\dots\Delta_{h,j}}\limits^{l}f(\boldsymbol{x}).
$$

Нехай задано вектори $\boldsymbol{h}=(h_1,\dots,h_d)$, $h_j \in
\mathbb{R}$, та ${\boldsymbol{k}=(k_1,\dots,k_d)}$, $k_j\in
\mathbb{Z}_+$, $j=\overline{1,d}$. Тоді мішана різниця
$\boldsymbol{k}$-го порядку з векторним кроком $\boldsymbol{h}$
визначається рівністю
$$
 \Delta_{\boldsymbol{h}}^{\boldsymbol{k}}f(\boldsymbol{x})=
 \Delta_{h_1,1}^{k_1}\Delta_{h_2,2}^{k_2}\dots\Delta_{h_d,d}^{k_d}f(\boldsymbol{x}).
$$
Крім цього покладемо $e_d=\{1,2, ... ,d\}$,  ${d\in \mathbb{N}}$, і
${e=\{j_1, ... ,j_m \}}$, ${m\in \mathbb{N}}$, ${m\leqslant d}$,
${1\leqslant j_1<j_2< ... <j_m\leqslant d}$. Задамо
невід'ємний вектор ${\boldsymbol{r}^e=(r_{j_1},\dots,r_{j_m})}$,
${r_j\geqslant 0}$, ${j=\overline{1,d}}$, і
${\bar{\boldsymbol{r}}^e=(\bar{r}_1,\dots,\bar{r}_d)}$, де
$$
 \bar{r}_i=
 \begin{cases}
    r_i, & i \in e; \\
    0, & i \in e_d\backslash e.
 \end{cases}
$$

Простори $S_{p,\theta}^{\boldsymbol{r}}B(\mathbb{R}^d)$, $1\leqslant
p,\theta \leqslant \infty$, де $\boldsymbol{r}$~--- заданий  вектор
із невід'ємними координатами означаються таким чином~\cite{Amanov_1965}:

1) якщо $1\leqslant \theta < \infty$, то
$$
 S_{p,\theta}^{\boldsymbol{r}}B(\mathbb{R}^d)=\Big \{f \in
 L_p(\mathbb{R}^d): \|f\|_{S_{p,\theta}^{\boldsymbol{r}} B(\mathbb{R}^d)} <\infty \Big\},
$$
де норма задається рівністю
$$
 \|f\|_{S_{p,\theta}^{\boldsymbol{r}} B(\mathbb{R}^d)}=\|f\|_p+
 \sum \limits_{ e \subset e_d \atop e \neq \varnothing} \left( \int \limits_0^2 \dots \int \limits_0^2 \prod \limits_{j\in e}
 h_j^{-\theta r_j-1}\|\Delta_{\boldsymbol{h}^e}^{\boldsymbol{k}^e}f(\cdot)\|_p^{\theta}\prod \limits_{j\in e}
 dh_j\right)^{\frac{1}{\theta}};
$$

2) якщо $\theta=\infty$, то
$$
S_{p,\infty}^{\boldsymbol{r}}B(\mathbb{R}^d)= \Big\{f \in
L_{p}(\mathbb{R}^d):
 \|f\|_{S_{p,\infty}^{\boldsymbol{r}}B(\mathbb{R}^d)} < \infty \Big\}
$$
i
$$
 \|f\|_{S_{p,\infty}^{\boldsymbol{r}}B(\mathbb{R}^d)}= \|f\|_p+
 \sum \limits_{ e \subset e_d \atop e \neq \varnothing} \sup \limits_{\boldsymbol{h}>0}\prod \limits_{j\in e}
 h_j^{-r_j}\|\Delta_{\boldsymbol{h}^e}^{\boldsymbol{k}^e}f(\cdot)\|_p,
$$
де $k_j>r_j\geqslant 0$, $j=\overline{1,d}$. Зазначимо, що простори функцій
$S_{p,\theta}^{\boldsymbol{r}}B(\mathbb{R}^d)$ при значенні параметра
$\theta=\infty$ збігаються з просторами
$S^{\boldsymbol{r}}_{p}H(\mathbb{R}^d)$, які вперше розглянув
С.\,М.~Нікольський~\cite{Nikolsky_63}, а у випадку $1\leqslant \theta< \infty$ вони були введені Т.\,І.~Амановим~\cite{Amanov_1965}.

Далі, замість
$S_{p,\theta}^{\boldsymbol{r}}B(\mathbb{R}^d)$ і
$S^{\boldsymbol{r}}_{p}H(\mathbb{R}^d)$  часто будемо використовувати  позначення
$S_{p,\theta}^{\boldsymbol{r}}B$ та $S^{\boldsymbol{r}}_{p}H$
відповідно.

У подальшому будемо вважати, що координати вектора
$\boldsymbol{r}=(r_1,\dots,r_d)$ впорядковані таким чином:
$0<r_1=r_2=\dots=r_{\nu}<r_{\nu+1}\leqslant\dots\leqslant r_d$.
Вектору $\boldsymbol{r}=(r_1,\dots,r_d)$ поставимо у відповідність вектор
$\boldsymbol{\gamma}=(\gamma_1,\dots,\gamma_d)$, $\gamma_j=r_j/r_1$, $j=\overline{1,d}$, а вектору $\boldsymbol{\gamma}$~--- $\boldsymbol{\gamma}'$, де
$\gamma'_j=\gamma_j$, якщо $j=\overline{1,\nu}$ i
$1<\gamma_j'<\gamma_j$, $j=\overline{\nu+1,d}$.

Також дамо означення просторів Нікольського--Бєсова функцій мішаної гладкості $S^{\boldsymbol{r}}_{p,\theta}B(\mathbb{R}^d)$ опосередковано через, так зване, декомпозиційне представлення елементів цих просторів. Уперше декомпозиційне представлення та відповідне йому нормування з'явилося у роботі   С.\,М.~Нікольського та П.\,І.~Лізоркіна~\cite{Lizorkin_Nikolsky_1989} і, як з'ясувалося пізніше, відіграло ключову роль у дослідженнях, які пов'язані з апроксимацією класів функцій. Наведемо спочатку необхідні означення та позначення.

Нехай $S=S(\mathbb{R}^d)$~--- простір Л.~Шварца основних нескінченно
диференційовних на $\mathbb{R}^d$ комплекснозначних функцій
$\varphi$, що спадають на нескінченності разом зі своїми похідними
швидше за будь-який степінь функції $\left(x_1^2+\ldots+x_d^2\right)^{-\frac{1}{2}}$ (див., наприклад, \cite{Lizorkin_69}). Через $S'$ позначимо
простір лінійних неперервних функціоналів над $S$. Зазначимо, що
елементами простору $S'$ є узагальнені функції.

Через  $\mathfrak{F}\varphi$ та $\mathfrak{F}^{-1}\varphi$ будемо позначати відповідно пряме та обернене перетворення Фур'є функцій з просторів $S$ та $S'$.

Носієм узагальненої функції $f$ будемо називати замикання
$\overline{\mathfrak{N}}$ такої множини точок
$\mathfrak{N}\subset\mathbb{R}^d$, що для довільної $\varphi \in S$,
яка дорівнює нулю в $\overline{\mathfrak{N}}$, виконується рівність
$\langle f,\varphi \rangle = 0$. Носій узагальненої функції $f$
будемо позначати через $\mbox{supp}\, f$.

Зазначимо, що для $1\leqslant p \leqslant \infty$ існує природне неперервне
вкладення $L_p(\mathbb{R}^d)$ в $S'$ і в цьому сенсі функції з
$L_p(\mathbb{R}^d)$ ототожнюються з елементами з $S'$.

Далі нехай $K_m(t)=\int_{\mathbb{R}} k_m(\lambda)e^{-2\pi i \lambda t} d\lambda$, $m\in \mathbb{Z}_+$, $K_{-1}\equiv 0$,  де

\begin{minipage}{9 cm}
$$
 k_m(\lambda)=
 \begin{cases}
    1, & |\lambda|<2^{m-1}, \\
    2(1-\frac{|\lambda|}{2^m}), & 2^{m-1}\leqslant|\lambda|\leqslant2^{m}, \\
    0, &  |\lambda|>2^{m},
 \end{cases}
$$
\end{minipage}
\begin{minipage}{7 cm}
$$
 k_0(\lambda)=
 \begin{cases}
    1-|\lambda|, & 0\leqslant|\lambda|\leqslant1, \\
    0, & |\lambda|>1.
 \end{cases}
$$
\end{minipage}

Для кожного вектора ${\boldsymbol{s}=(s_{1},...,s_{d})}$, ${s_{j}\in
\mathbb{Z}_+}$, ${j=\overline{1,d}}$, покладемо
\begin{equation}\label{A_s_prod}
A^*_{\boldsymbol{s}}(\boldsymbol{x})=\prod\limits_{j=1}^d \big(K_{s_j}(x_j)-K_{s_j-1}(x_j)\big),
\end{equation}
$$
A^*_{\boldsymbol{s}}(f,\boldsymbol{x})=f(\boldsymbol{x})\ast A^*_{\boldsymbol{s}}(\boldsymbol{x})=
\int\limits_{\mathbb{R}^d}f(\boldsymbol{y}) A^*_{\boldsymbol{s}}(\boldsymbol{x}-\boldsymbol{y}) d\boldsymbol{y}.
$$
Одразу зауважимо, що $A^*_{\boldsymbol{s}}(f,\boldsymbol{x})$~--- блоки Валле Пуссена функції $f$.
Також для $\boldsymbol{s}\in \mathbb{Z}^d_+$ розглянемо множини
$$
Q_{2^{\boldsymbol{s}}}^*=\big\{\boldsymbol{\lambda}=(\lambda_1,...,\lambda_d):
\ \eta(s_j)2^{s_{j}-1}\leqslant |\lambda_j|<2^{s_j}, \lambda_j\in
\mathbb{R}, \ \ j=\overline{1,d}\big\},
$$
$$
\rho_+(\boldsymbol{s}):=\big\{\boldsymbol{k}=(k_1,...,k_d):
\eta(s_j)2^{s_j-1}\leqslant k_j<2^{s_j}, k_{j}\in\mathbb{Z}_+,
j=\overline{1,d}\big\},
$$
де $\eta(0)=0$ і $\eta(t)=1, \ t>0$.

Має місце таке твердження.

\bf Лема А \rm  (див., наприклад, \cite{WangHeping_SunYongsheng_1999_AppT})\textbf{.} \it Нехай $1\leqslant p\leqslant\infty$, тоді для будь-якої функції $f\in
L_p(\mathbb{R}^d)$ маємо
 $$
f(\boldsymbol{x})\stackrel{L_p}{=}\sum A^*_{\boldsymbol{s}}(f,\boldsymbol{x})
 $$
 і крім того \rm $\mbox{supp} \, \mathfrak{F}A_{\boldsymbol{s}}(f,\boldsymbol{x})\subseteq Q_{2^{\boldsymbol{s}}}^*$. \rm

 У прийнятих позначеннях простори $S_{p,\theta}^{\boldsymbol{r}} B(\mathbb{R}^d)$,
$1\leqslant p, \theta \leqslant \infty$, $\boldsymbol{r}>0$, можна означити таким
чином~(див., наприклад, \cite{WangHeping_SunYongsheng_1999_AppT}, \cite{WangHeping_1997_Q}):
$$
S^{\boldsymbol{r}}_{p,\theta}B:=\Big\{f\in L_p(\mathbb{R}^d): \
\|f\|_{S^{\boldsymbol{r}}_{p,\theta}B}<\infty \Big\},
$$
де
\begin{equation}\label{Norm_dek_Sr1}
   \|f\|_{S^{\boldsymbol{r}}_{p,\theta}B}\asymp \Bigg(\sum \limits_{\boldsymbol{s}\geqslant
   0}2^{(\boldsymbol{s},\boldsymbol{r})\theta}
   \|A^*_{\boldsymbol{s}}(f,\cdot)\|_p^{\theta}\Bigg)^{\frac{1}{\theta}}
\end{equation}
  при $1\leqslant\theta<\infty$ і
\begin{equation}\label{Norm_dek_Sr1_inf}
   \|f\|_{S^{\boldsymbol{r}}_{p}H}\asymp \sup
   \limits_{\boldsymbol{s}\geqslant
   0} 2^{(\boldsymbol{s},\boldsymbol{r})}\|A^*_{\boldsymbol{s}}(f,\cdot)\|_p.
\end{equation}
\rm

Тут і надалі по тексту для додатних величин $A$ і  $B$ вживається запис  $A\asymp B$, який означає, що існують такі додатні  сталі $C_1$ та $C_2$, які не залежать від одного істотного параметра у величинах  $A$ і  $B$ (наприклад, у вище наведених співвідношеннях (\ref{Norm_dek_Sr1}) і (\ref{Norm_dek_Sr1_inf})~--- від функції $f$), що ${C_1 A \leqslant B \leqslant C_2 A}$. Якщо тільки $B\leqslant C_2 A $ $\big(B \geqslant C_1 A\big)$, то пишемо
$B\ll A$ $\big(B \gg A \big)$. Всі сталі $C_i$, $i=1,2,...$, які
зустрічаються у роботі, залежать, можливо, лише від параметрів, що
входять в означення класу, метрики, в якій оцінюється похибка
наближення, та розмірності простору $\mathbb{R}^d$.

Окрім цього нагадаємо, що у випадку $1<p<\infty$ норму функцій з просторів $S_{p,\theta}^{\boldsymbol{r}} B(\mathbb{R}^d)$ можна означити в дещо іншій формі.

Нехай $A\subset \mathbb{R}^d$~--- деяка вимірна множина. Позначимо через
$\chi_{_A}$ характеристичну функцію множини $A$ і для ${f \in
L_p(\mathbb{R}^d)}$ покладемо
$$
\delta_{\boldsymbol{s}}^*(f,\boldsymbol{x})=\mathfrak{F}^{-1}(\chi_{Q_{2^{\boldsymbol{s}}}^*}\cdot
\mathfrak{F}f).
$$

Тоді простори $S_{p,\theta}^{\boldsymbol{r}} B$,
$1<p<\infty$, $1\leqslant  \theta \leqslant \infty$, $\boldsymbol{r}>0$, можна означити таким
чином~\cite{Lizorkin_Nikolsky_1989}:
$$
S^{\boldsymbol{r}}_{p,\theta}B:=\Big\{f\in L_p(\mathbb{R}^d): \
\|f\|_{S^{\boldsymbol{r}}_{p,\theta}B}<\infty \Big\},
$$
де
\begin{equation}\label{Norm_dek_Sr}
   \|f\|_{S^{\boldsymbol{r}}_{p,\theta}B}\asymp \Bigg(\sum \limits_{\boldsymbol{s}\geqslant
   0}2^{(\boldsymbol{s},\boldsymbol{r})\theta}
   \|\delta_{\boldsymbol{s}}^*(f,\cdot)\|_p^{\theta}\Bigg)^{\frac{1}{\theta}}
\end{equation}
  при $1\leqslant\theta<\infty$ і
\begin{equation}\label{Norm_dek_Sr_inf}
   \|f\|_{S^{\boldsymbol{r}}_{p}H}\asymp \sup
   \limits_{\boldsymbol{s}\geqslant
   0} 2^{(\boldsymbol{s},\boldsymbol{r})}\|\delta_{\boldsymbol{s}}^*(f,\cdot)\|_p.
\end{equation}
\rm

Під класом $S^{\boldsymbol{r}}_{p,\theta}B$ будемо розуміти множину
функцій ${f \in L_p(\mathbb{R}^d)}$ для яких
$\|f\|_{S^{\boldsymbol{r}}_{p,\theta}B}\leqslant  1$ і при цьому збережемо
для класів $S^{\boldsymbol{r}}_{p,\theta}B$ ті ж самі позначення, що і для
просторів $S^{\boldsymbol{r}}_{p,\theta}B$.

Як видно з (\ref{Norm_dek_Sr1})\,--\,(\ref{Norm_dek_Sr_inf}), для $f\in S^{\boldsymbol{r}}_{p,\theta}B$ при деякому значенні $p$, $1<p<\infty$, має місце співвідношення
\begin{equation}\label{As_deltas}
 \|\delta_{\boldsymbol{s}}^*(f,\cdot)\|_p\asymp \|A^*_{\boldsymbol{s}}(f,\cdot)\|_p.
\end{equation}

Перейдемо до означення апроксимативних характеристик, які розглядаються у роботі.

Нехай $\mathcal{L}\subset \mathbb{Z}^d_+$~--- деяка скінченна
множина. Покладемо
$$
Q(\mathcal{L})=\bigcup\limits_{\boldsymbol{s}\in \mathcal{L}}Q_{2^{\boldsymbol{s}}}^*
$$
і позначимо
$$
G\big(Q(\mathcal{L})\big)=\Big\{f\in L_q(\mathbb{R}^d): \
\mbox{supp}\mathfrak{F}f\subseteq Q(\mathcal{L}) \Big\}.
$$
Відомо, що елементами множини $G\big(Q(\mathcal{L})\big)$ є цілі функції експоненціального типу.

Для $f\in L_q(\mathbb{R}^d)$, $1\leqslant q \leqslant \infty$,
означимо величину
$$
E\big(f, G\big(Q(\mathcal{L})\big)\big)_q:=E_{Q(\mathcal{L})}(f)_q:=
\inf\limits_{g\in  G(Q(\mathcal{L}))}\|f(\cdot)-g(\cdot)\|_q,
$$
яка називається найкращим наближенням функції $f$ цілими функціями з
множини $G\big(Q(\mathcal{L})\big)$. Якщо $F\subset
L_q(\mathbb{R}^d)$~--- деякий функціональний клас, то покладемо
\begin{equation}\label{EQN}
E_{Q(\mathcal{L})}(F)_q=\sup\limits_{f\in F}E_{Q(\mathcal{L})}(f)_q.
\end{equation}

Далі для $f\in L_q(\mathbb{R}^d)$, $1\leqslant q \leqslant \infty$, покладемо
$$
S_{Q(\mathcal{L})}f(\boldsymbol{x})= S_{Q(\mathcal{L})}(f,\boldsymbol{x})=\sum\limits_{\boldsymbol{s}\in
\mathcal{L}}\delta_{\boldsymbol{s}}^*(f,\boldsymbol{x}), \ \ \boldsymbol{x}\in  \mathbb{R}^d
$$
і означимо
\begin{equation}\label{EQN1}
\mathcal{E}_{Q(\mathcal{L})}(f)_q=
\|f(\cdot)-S_{Q(\mathcal{L})}f(\cdot)\|_q \ \ \mbox{та} \ \ \mathcal{E}_{Q(\mathcal{L})}(F)_q=\sup\limits_{f\in F}\mathcal{E}_{Q(\mathcal{L})}(f)_q.
\end{equation}

Наші дослідження величин (\ref{EQN}) та (\ref{EQN1}) проводяться у випадку, коли ${F=S^{\boldsymbol{r}}_{p,\theta}B(\mathbb{R}^d)}$, а множина $Q(\mathcal{L})$ визначається таким чином:
$$
Q(\mathcal{L})=Q_n^{\boldsymbol{\gamma}}=\bigcup \limits_{(\boldsymbol{s},\boldsymbol{\gamma}) \leqslant n}Q_{2^{\boldsymbol{s}}}^*,
$$
де $n\in \mathbb{N}$. Множина $Q_n^{\boldsymbol{\gamma}}$ породжує в
$\mathbb{R}^d$, так званий, східчастий гіперболічний хрест.

Насамперед зазначимо, що при $1< q < \infty$ і  $f\in L_q(\mathbb{R}^d)$
 має  місце
співвідношення (див., наприклад,~\cite{Lizorkin_Nikolsky_1989})
\begin{equation}\label{e=e}
E_{Q(\mathcal{L})}(f)_q\leqslant\mathcal{E}_{Q(\mathcal{L})}(f)_q\leqslant C_3 E_{Q(\mathcal{L})}(f)_q,
\end{equation}
де $C_3 \geqslant1$~---  деяка стала.

\newpage

\vskip 3 mm  \textbf{2. Допоміжні твердження}.

\bf Теорема А \rm \cite{Amanov_1965}. \it Нехай $1\leqslant  p, \theta \leqslant \infty$,
$1\leqslant  p\leqslant  q \leqslant \infty$ і маємо такий вектор $\boldsymbol{\rho}$, що
$\rho_j=r_j-\left(\frac{1}{p}-\frac{1}{q}\right)>0$, $j=\overline{1,d}$. Тоді, якщо $f\in S_{p,\theta}^{\boldsymbol{r}}B(\mathbb{R}^d)$, то $f\in   S_{q,\theta}^{\boldsymbol{\rho}}B(\mathbb{R}^d)$
 і
 $$
  \|f\|_{S_{q,\theta}^{\boldsymbol{\rho}}B(\mathbb{R}^d)}\ll \|f\|_{S_{p,\theta}^{\boldsymbol{r}}B(\mathbb{R}^d)}.
 $$

\vskip 1 mm

\bf Теорема~Б \rm \cite[c.~150]{Nikolsky_1969_book} \textbf{.} \it
Якщо $1\leqslant p \leqslant q \leqslant \infty$, тоді для цілої
функції експоненціального типу $g_{\boldsymbol{\nu}}\in L_p(\mathbb{R}^d)$, $\boldsymbol{\nu}=(\nu_1,\ldots,\nu_2)$, $\nu_i\geqslant 0$, $i=\overline{1,d}$, має
місце нерівність (різних метрик)
$$
 \|g_{\boldsymbol{\nu}}\|_{L_{q}(\mathbb{R}^d)}\leqslant 2^d\left( \prod \limits_{j=1}^d
 \nu_k\right)^{\frac{1}{p}-\frac{1}{q}}\|g_{\boldsymbol{\nu}}\|_{L_p(\mathbb{R}^d)}.
$$\rm

\vskip 1 mm \textbf{Теорема~В} \rm (Літтлвуда\,--\,Пелі) (див.,
наприклад, \cite[c.~81]{Nikolsky_1969_book})\textbf{.} \it Нехай задано
${1\!<p<\infty}$. Існують такі додатні числа $C_{4}, C_{5}$,  що для
кожної функції ${f\in L_{p}(\mathbb{R}^d)}$ виконуються
співвідношення
 $$
 C_{4}\|f\|_{p}\leqslant  \Bigg\| \Bigg(\sum\limits_{\boldsymbol{s}\geqslant 0}
  |\delta_{\boldsymbol{s}}^*(f,\cdot)|^{2} \Bigg)^{\frac{1}{2}} \Bigg\|_{p}\leqslant  C_{5}\|f\|_{p}.
$$ \rm

\vskip 1mm

\bf Лема Б  \rm \cite{WangHeping_SunYongsheng_1995}\textbf{.} \it Нехай задано $1<p<q<\infty$ і $f\in
L_q(\mathbb{R}^d)$. Тоді
 $$
\|f\|_q\ll \bigg( \sum \limits_{\boldsymbol{s}\geqslant 0} \|\delta_{\boldsymbol{s}}^*(f,\cdot)\|_p^q
\; 2^{\|\boldsymbol{s}\|_1(\frac{1}{p}-\frac{1}{q})q}\bigg)^{\frac{1}{q}},
 $$
 де $\|\boldsymbol{s}\|_1=s_1+\ldots + s_d$.
 \rm

Лема~Б є аналогом  леми, яка вперше була доведена для періодичного випадку
\cite[c.~25]{Temlyakov_1986m}.

\bf Лема В \rm \cite[с.~11]{Temlyakov_1986m}\textbf{.} \it Має місце оцінка
 $$
 \sum\limits_{(\boldsymbol{s},\boldsymbol{\gamma})\geqslant n}2^{-\alpha(\boldsymbol{s},\boldsymbol{\gamma})}\asymp 2^{-\alpha
 n}n^{d-1}, \ \ \alpha>0
 $$

\vskip 1 mm

 \bf Лема Г \rm \cite[с.~11]{Temlyakov_1986m}\textbf{.} \it
Має місце оцінка
 $$
 \sum\limits_{(\boldsymbol{s},\boldsymbol{\gamma'})\geqslant n}2^{-\alpha(\boldsymbol{s},\boldsymbol{\gamma})}\asymp 2^{-\alpha
 n}n^{\nu-1}, \ \ \alpha>0.
 $$
\rm \vskip 2.5 mm

\textbf{3. Наближення цілими функціями з носіями їх перетворення Фур'є у східчастому гіперболічному хресті.} Перш ніж перейти до формулювання та доведення основних результатів, встановимо декілька допоміжних оцінок.

\vskip 1 mm

 \bf Лема 1. \it
Для $A_{\boldsymbol{s}}(\boldsymbol{x})$ мають місце такі оцінки
 \begin{equation}\label{As_infty}
 \|A^*_{\boldsymbol{s}}(\cdot)\|_{\infty}\asymp 2^{\|\boldsymbol{s}\|_1},
 \end{equation}
\begin{equation}\label{As_Sum_infty}
 \bigg\|\sum\limits_{(\boldsymbol{s},1)=n+1}A^*_{\boldsymbol{s}}(\cdot)\bigg\|_{\infty}\asymp 2^n n^{d-1}.
 \end{equation}
\vskip 1 mm \rm

{\textbf{\textit{Доведення.}}} \ \rm Врахувавши, що $A^*_{\boldsymbol{s}}(\boldsymbol{x})$ визначається згідно з формулою (\ref{A_s_prod}) та провівши деякі перетворення для $K_{s_j}(x_j)$, отримаємо
$$
A^*_{\boldsymbol{s}}(\boldsymbol{x})=\prod\limits_{j=1}^d \frac{2^{2-s_j}\sin^2\pi 2^{s_j-2}x_j \big(2\cos\pi 2^{s_j-1}x_j+1\big) }{\pi^2 x_j^2} \times
$$
$$
\times\big(\cos\pi 2^{s_j-1}x_j+\cos\pi 2^{s_j}x_j-1\big)=\prod\limits_{j=1}^d \mathfrak{I}_j(x_j).
$$

Тоді для норми $A^*_{\boldsymbol{s}}(\boldsymbol{x})$ у просторі $L_{\infty}(\mathbb{R}^d)$ можемо записати
\begin{equation}\label{As_infty_1}
\|A^*_{\boldsymbol{s}}(\cdot)\|_{\infty}= \sup\limits_{\boldsymbol{x}\in \mathbb{R}^d}  \bigg|\prod\limits_{j=1}^d \mathfrak{I}_j(x_j)\bigg| =\prod\limits_{j=1}^d \sup\limits_{x_j\in \mathbb{R}} |\mathfrak{I}_j(x_j)|.
\end{equation}

Таким чином для оцінки норми $A^*_{\boldsymbol{s}}(\boldsymbol{x})$ достатньо буде оцінити $\sup\limits_{x_j\in \mathbb{R}} |\mathfrak{I}_j(x_j)|$, $j=\overline{1,d}$.

Встановимо спочатку оцінку зверху. Зробимо заміну $\pi 2^{s_j-2}x_j=t_j$, отримаємо
$$
\sup\limits_{x_j\in \mathbb{R}} |\mathfrak{I}_j(x_j)|=\sup\limits_{t_j\in \mathbb{R}}\left|\frac{2^{s_j-2}\sin^2 t_j \big(2\cos 2t_j+1\big) \big(\cos 2t_j+\cos 4t_j-1\big)}{t_j^2}\right| \ll
$$
$$
\ll 2^{s_j}\sup\limits_{t_j\in \mathbb{R}}\left| 9 \, \frac{\sin^2 t_j }{t_j^2}\right| \ll 2^{s_j}.
$$

Скориставшись останньою оцінкою і (\ref{As_infty_1}), одержуємо
$$
\|A^*_{\boldsymbol{s}}(\cdot)\|_{\infty} \ll \prod\limits_{j=1}^d 2^{s_j}= 2^{\|\boldsymbol{s}\|_1}.
$$

Для оцінки знизу маємо
$$
\sup\limits_{x_j\in \mathbb{R}} |\mathfrak{I}_j(x_j)|=\sup\limits_{t_j\in \mathbb{R}}\left|\frac{2^{s_j-2}\sin^2 t_j \big(2\cos 2t_j+1\big) \big(\cos 2t_j+\cos 4t_j-1\big)}{t_j^2}\right| \gg
$$
\begin{equation}\label{As_infty_lower}
\gg 2^{s_j} \left|\frac{\sin^2 \frac{\pi}{12} \big(2\cos \frac{\pi}{6}+1\big) \big(\cos \frac{\pi}{6}+\cos \frac{\pi}{3}-1\big)}{\left(\frac{\pi}{12}\right)^2}\right| \gg 2^{s_j}.
\end{equation}

Підставивши (\ref{As_infty_lower}) в (\ref{As_infty_1}), одержуємо
$$
\|A^*_{\boldsymbol{s}}(\cdot)\|_{\infty} \gg \prod\limits_{j=1}^d 2^{s_j}= 2^{\|\boldsymbol{s}\|_1}.
$$

Отже, оцінку $(\ref{As_infty})$ встановлено.

Перейдемо до встановлення оцінки $(\ref{As_Sum_infty})$. Оцінка зверху безпосередньо випливає з $(\ref{As_infty})$ та нерівності Мінковського. Дійсно
$$
\bigg\|\sum\limits_{(\boldsymbol{s},1)=n+1}A^*_{\boldsymbol{s}}(\cdot)\bigg\|_{\infty}\ll
\sum\limits_{(\boldsymbol{s},1)=n+1} \|A^*_{\boldsymbol{s}}(\cdot)\|_{\infty}\ll \sum\limits_{(\boldsymbol{s},1)=n+1} 2^{\|\boldsymbol{s}\|_1} \asymp 2^n n^{d-1}.
$$

Для оцінки знизу маємо
$$
\bigg\|\sum\limits_{(\boldsymbol{s},1)=n+1}A^*_{\boldsymbol{s}}(\cdot)\bigg\|_{\infty}= \sup\limits_{\boldsymbol{x}\in \mathbb{R}^d} \bigg|\sum\limits_{(\boldsymbol{s},1)=n+1}A^*_{\boldsymbol{s}}(\boldsymbol{x})\bigg|=
$$
$$
=\sup\limits_{\boldsymbol{t}\in \mathbb{R}^d} \left|\sum\limits_{(\boldsymbol{s},1)=n+1} \prod\limits_{j=1}^d \frac{2^{s_j-2}\sin^2 t_j \big(2\cos 2t_j+1\big) \big(\cos 2t_j+\cos 4t_j-1\big)}{t_j^2}\right|\gg
$$
$$
\gg \sum\limits_{(\boldsymbol{s},1)=n+1} \prod\limits_{j=1}^d  2^{s_j} \left|\frac{\sin^2 \frac{\pi}{12} \big(2\cos \frac{\pi}{6}+1\big) \big(\cos \frac{\pi}{6}+\cos \frac{\pi}{3}-1\big)}{\left(\frac{\pi}{12}\right)^2}\right| \gg
$$
$$
\gg \sum\limits_{(\boldsymbol{s},1)=n+1} 2^{\|\boldsymbol{s}\|_1} \asymp 2^n n^{d-1}.
$$

Оцінку $(\ref{As_Sum_infty})$ встановлено. Лему~1 доведено.

\vskip 1 mm
 \bf Лема 2. \it Нехай $1\leqslant p<\infty$, тоді має місце оцінка
 \begin{equation}\label{As_p}
 \|A^*_{\boldsymbol{s}}(\cdot)\|_{p}\asymp 2^{\|\boldsymbol{s}\|_1\left(1-\frac{1}{p}\right)}.
 \end{equation}

\vskip 1 mm \rm

{\textbf{\textit{Доведення.}}} \ \rm Згідно з позначеннями леми~1, можемо записати
$$
\|A^*_{\boldsymbol{s}}(\cdot)\|_{p}= \bigg\|\prod\limits_{j=1}^d \mathfrak{I}_j(\cdot)\bigg\|_p =\left(\int\limits_{\mathbb{R}^d}\bigg|\prod\limits_{j=1}^d \mathfrak{I}_j(x_j)\bigg|^p d\boldsymbol{x}\right)^{\frac{1}{p}}=
$$
$$
=\left(\int\limits_{\mathbb{R}^d}\prod\limits_{j=1}^d |\mathfrak{I}_j(x_j)|^p d\boldsymbol{x}\right)^{\frac{1}{p}}=\left(\prod\limits_{j=1}^d\int\limits_{\mathbb{R}} |\mathfrak{I}_j(x_j)|^p dx_j\right)^{\frac{1}{p}}=
$$
$$
=\left(\prod\limits_{j=1}^d\int\limits_{\mathbb{R}} \Big|\frac{2^{2-s_j}\sin^2\pi 2^{s_j-2}x_j \big(2\cos\pi 2^{s_j-1}x_j+1\big)}{\pi^2 x_j^2} \times \right.
$$
\begin{equation}\label{As_infty_2}
 \times \big(\cos\pi 2^{s_j-1}x_j+\cos\pi 2^{s_j}x_j-1\big)\Big|^p dx_j\Bigg)^{\frac{1}{p}}
\end{equation}

Оцінимо зверху інтеграл у співвідношенні (\ref{As_infty_2}). Зробивши заміну ${\pi 2^{s_j-2}x_j=t_j}$, отримаємо
$$
\int\limits_{\mathbb{R}} \Big|\frac{2^{2-s_j}\sin^2\pi 2^{s_j-2}x_j \big(2\cos\pi 2^{s_j-1}x_j+1\big) \big(\cos\pi 2^{s_j-1}x_j+\cos\pi 2^{s_j}x_j-1\big)}{\pi^2 x_j^2}\Big|^p dx_j \!\!\leqslant
$$
$$
\leqslant \int\limits_{\mathbb{R}} \Big|\frac{2^{2-s_j}\sin^2\pi 2^{s_j-2}x_j }{\pi^2 x_j^2}\Big|^p dx_j = \int\limits_{\mathbb{R}}\frac{2^{(s_j-2)(p-1)}}{\pi} \Big|\frac{\sin^2t_j }{t_j^2}\Big|^p dt_j
$$
\begin{equation}\label{As_int}
\ll 2^{s_j(p-1)}.
\end{equation}

Підставивши (\ref{As_int}) в (\ref{As_infty_2}), одержуємо
$$
\|A^*_{\boldsymbol{s}}(\cdot)\|_p \ll \left(\prod\limits_{j=1}^d 2^{s_j(p-1)}\right)^{\frac{1}{p}}= 2^{\|\boldsymbol{s}\|_1\left(1-\frac{1}{p}\right)}.
$$

Тепер оцінимо знизу інтеграл у співвідношенні (\ref{As_infty_2}). Зробивши заміну ${\pi 2^{s_j-2}x_j=t_j}$, отримаємо
$$
\int\limits_{\mathbb{R}} \Big|\frac{2^{2-s_j}\sin^2\pi 2^{s_j-2}x_j \big(2\cos\pi 2^{s_j-1}x_j+1\big) \big(\cos\pi 2^{s_j-1}x_j+\cos\pi 2^{s_j}x_j-1\big)}{\pi^2 x_j^2}\Big|^p dx_j\!\!=
$$
$$
=\int\limits_{\mathbb{R}} 2^{(s_j-2)(p-1)}\Big|\frac{\sin^2 t_j \big(2\cos 2t_j+1\big) \big(\cos 2t_j+\cos 4t_j-1\big)}{t_j^2}\Big|^p dt_j \gg
$$
$$
\gg2^{(s_j-2)(p-1)} \int\limits_{\frac{7\pi}{12}}^{\frac{5\pi}{12}} \Big|\frac{\sin^2 t_j \big(2\cos 2t_j+1\big) \big(\cos 2t_j+\cos 4t_j-1\big)}{t_j^2}\Big|^p dt_j \gg 2^{s_j(p-1)}.
$$
Врахувавши, що підінтегральна функція неперервна на проміжку $\big[\frac{5\pi}{12},\frac{7\pi}{12}\big]$ та досягає свого найменшого значення, одержуємо
\begin{equation}\label{As_int2}
\gg 2^{s_j(p-1)}\int\limits_{\frac{\pi}{3}}^{\frac{2\pi}{3}} M_j^p dt_j\asymp 2^{s_j(p-1)}.
\end{equation}

Скориставшись (\ref{As_int2}), одержуємо
$$
\|A^*_{\boldsymbol{s}}(\cdot)\|_p \gg \left(\prod\limits_{j=1}^d 2^{s_j(p-1)}\right)^{\frac{1}{p}}= 2^{\|\boldsymbol{s}\|_1\left(1-\frac{1}{p}\right)}.
$$

Отже, оцінку $(\ref{As_p})$ встановлено. Лему~2 доведено.

\vskip 3 mm  \bf Теорема~1. \it Нехай $r_1>1$, $1 \leqslant \theta \leqslant \infty$. Тоді має
місце порядкове співвідношення
\begin{equation} \label{teor_1}
  \mathcal{E}_{Q_n^{\gamma}}\big(S^{\boldsymbol{r}}_{1,\theta}B\big)_{\infty}= \sup \limits_{f\in
  S_{1,\theta}^{\boldsymbol{r}}B}\|f(\cdot)-S_{Q_n^{\boldsymbol{\gamma}}}(f,\cdot)\|_{\infty}\asymp
  2^{-n\left(r_1-1\right)}n^{(\nu-1)\left(1-\frac{1}{\theta}\right)}.
\end{equation}
\vskip 1 mm \rm

{\textbf{\textit{Доведення.}}} \ \rm  Встановимо спочатку оцінку
зверху в (\ref{teor_1}). Нехай $f\in S^{\boldsymbol{r}}_{1,\theta}B$. Оскільки $r_1>1$,  то на основі теорема~А можемо стверджувати, що при деякому $1<q_0<\infty$, $f\in S^{\boldsymbol{\rho}}_{q_0,\theta}B$, де $\rho_j=r_j-\left(1-\frac{1}{q_0}\right)>0$, $j=\overline{1,d}$.  Тоді, скориставшись нерівністю  Мінковського,  нерівністю різних метрик (теорема~Б) та (\ref{As_deltas}), можемо записати
$$
 \|f(\cdot)-S_{Q_n^{\boldsymbol{\gamma}}}(f,\cdot)\|_{\infty}=\Big\|f(\cdot)-
 \sum \limits_{(\boldsymbol{s},\boldsymbol{\gamma})\leqslant n}\delta_{\boldsymbol{s}}^*(f,\cdot)\Big\|_{\infty}\leqslant \sum
 \limits_{(\boldsymbol{s},\boldsymbol{\gamma})>n}\|\delta_{\boldsymbol{s}}^*(f,\cdot)\|_{\infty}\ll
$$
$$
\ll \sum
 \limits_{(\boldsymbol{s},\boldsymbol{\gamma})>n}2^{\frac{\|\boldsymbol{s}\|_1}{q_0}}
\|\delta_{\boldsymbol{s}}^*(f,\cdot)\|_{q_0}\asymp \sum
 \limits_{(\boldsymbol{s},\boldsymbol{\gamma})>n}2^{\frac{\|\boldsymbol{s}\|_1}{q_0}}
\|A_{\boldsymbol{s}}^*(f,\cdot)\|_{q_0}\ll
$$
\begin{equation} \label{os1_teor_1}
 \ll \sum
 \limits_{(\boldsymbol{s},\boldsymbol{\gamma})>n}2^{\frac{\|\boldsymbol{s}\|_1}{q_0}} 2^{\|\boldsymbol{s}\|_1\left(1-\frac{1}{q_0}\right)}
\|A_{\boldsymbol{s}}^*(f,\cdot)\|_{1}= \sum
 \limits_{(\boldsymbol{s},\boldsymbol{\gamma})>n}2^{\|\boldsymbol{s}\|_1}
\|A_{\boldsymbol{s}}^*(f,\cdot)\|_{1}.
\end{equation}

Щоб продовжити оцінку (\ref{os1_teor_1}) розглянемо спочатку
випадок, коли ${1\leqslant \theta<\infty}$. Тоді, застосувавши
нерівність Гельдера, з відповідною модифікацією при $\theta = \infty$, будемо мати
$$
 \sum
 \limits_{(\boldsymbol{s},\boldsymbol{\gamma})>n}2^{\|\boldsymbol{s}\|_1}
\|A_{\boldsymbol{s}}^*(f,\cdot)\|_{1}\leqslant \left(\sum
 \limits_{(\boldsymbol{s},\boldsymbol{\gamma})>n}
2^{(\boldsymbol{s},\boldsymbol{r})\theta}
\|A_{\boldsymbol{s}}^*(f,\cdot)\|_{1}^{\theta}\right)^{\frac{1}{\theta}}\times
$$
$$
\times\left(\sum
 \limits_{(\boldsymbol{s},\boldsymbol{\gamma})>n}2^{-(\boldsymbol{s},\boldsymbol{r}-1)
\frac{\theta}{\theta-1}}\right)^{1-\frac{1}{\theta}}\ll \|f\|_{S^{\boldsymbol{r}}_{1,\theta}B}\left(\sum
 \limits_{(\boldsymbol{s},\boldsymbol{\gamma})>n}2^{-(\boldsymbol{s},\boldsymbol{r}-1)
\frac{\theta}{\theta-1}}\right)^{1-\frac{1}{\theta}}\leqslant
$$
\begin{equation} \label{os2_teor2_1}
 \leqslant \left(\sum
 \limits_{(\boldsymbol{s},\boldsymbol{\gamma})>n}
2^{-(\boldsymbol{s},\boldsymbol{\bar{\gamma}})(r_1-1)
\frac{\theta}{\theta-1}}\right)^{1-\frac{1}{\theta}}=J_1,
\end{equation}
де $\boldsymbol{\bar{\gamma}}=(\bar{\gamma}_1,\dots,\bar{\gamma}_d)$~--- вектор з
координатами $\bar{\gamma}_j=(r_j-1)/(r_1-1)$, $j=\overline{1,d}$, а $\boldsymbol{r}-1$ позначає вектор з координатами
$r_j-1$, $j=\overline{1,d}$. Якщо $j=\overline{1,\nu}$, то
$\bar{\gamma}_j=\gamma_j$
  і $1<\gamma_j\leqslant \gamma_j$, якщо
$j=\overline{\nu+1,d}$. Тому, скориставшись лемою~В,
отримаємо оцінку
\begin{equation} \label{os3_teor2_1}
 J_1\ll2^{-n(r_1-1)}n^{(\nu-1)(1-\frac{1}{\theta})}.
\end{equation}

Отже, співставивши  (\ref{os1_teor_1})--(\ref{os3_teor2_1}),
отримаємо оцінку
$$
\sup \limits_{f\in
 S^{\boldsymbol{r}}_{1,\theta}B}\|f(\cdot)-S_{Q_n^{\boldsymbol{\gamma}}}(f,\cdot)\|_{\infty}\ll
 2^{-n(r_1-1)}n^{(\nu-1)(1-\frac{1}{\theta})}.
$$

Нехай тепер $\theta=\infty$. Тоді, згідно з означенням класів $S^{\boldsymbol{r}}_{1,\theta}B$ маємо ${\|A_{\boldsymbol{s}}^*(f,\cdot)\|_p\ll2^{-(\boldsymbol{s},\boldsymbol{r})}}$, скориставшись лемою~В,  для
(\ref{os1_teor_1}) можемо записати
$$
 \sum
\limits_{(\boldsymbol{s},\boldsymbol{\gamma})>n}2^{\|\boldsymbol{s}\|_1}
\|A_{\boldsymbol{s}}^*(f,\cdot)\|_1\ll\sum
\limits_{(\boldsymbol{s},\boldsymbol{\gamma})>n}2^{-(\boldsymbol{s},\boldsymbol{r}-1)}=
$$
\begin{equation} \label{os4_teot2_1}
=\sum
\limits_{(\boldsymbol{s},\boldsymbol{\gamma}>n)}
2^{-(\boldsymbol{s},\boldsymbol{\bar{\gamma}})(r_1-1)}\ll2^{-n(r_1-1)}n^{\nu-1}.
\end{equation}
Об'єднавши (\ref{os3_teor2_1}) та (\ref{os4_teot2_1}), одержуємо оцінку зверху в (\ref{teor_1}).

Перейдемо до встановлення оцінки знизу, яку достатньо отримати для
випадку $\nu=d$. Розглянемо функції
$$
 f_1(x)=C_6 2^{-nr_1} n^{-\frac{d-1}{\theta}} \sum
 \limits_{(\boldsymbol{s},1)=n+1}A^*_{\boldsymbol{s}}(\boldsymbol{x}), \ \ C_6>0,
$$
якщо $1\leqslant \theta<\infty$, i
$$
 f_2(x)=C_7 2^{-nr_1} \sum
 \limits_{(\boldsymbol{s},1)=n+1}A^*_{\boldsymbol{s}}(\boldsymbol{x}), \ \ C_7>0,
$$
якщо $\theta=\infty$.

Переконаємося, що дані функції належать класам $S^{\boldsymbol{r}}_{1,\theta}B$ і
$S^{\boldsymbol{r}}_{1,\infty}B$ відповідно. Оскільки, згідно з  лемою~1, при $p=1$ має місце оцінка $\big\|A_{\boldsymbol{s}}^*(\cdot)\big\|_1\asymp
 1$,  то
$$
 \|f_1\|_{S^{\boldsymbol{r}}_{1,\theta}B}\asymp \left( \sum
 \limits_{(\boldsymbol{s},1)=n+1}2^{(\boldsymbol{s},\boldsymbol{r})\theta}
\|A_{\boldsymbol{s}}^*(f_1,\cdot)\|_1^{\theta}\right)^{\frac{1}{\theta}}\asymp
$$
$$
 \asymp 2^{-nr_1} n^{-\frac{d-1}{\theta}} \left( \sum
 \limits_{(\boldsymbol{s},1)=n+1}2^{(\boldsymbol{s},\boldsymbol{r})\theta}\| A_{\boldsymbol{s}}^*(\cdot)\|_1^{\theta}\right)^{\frac{1}{\theta}}\asymp
$$
$$
 \asymp 2^{-nr_1} n^{-\frac{d-1}{\theta}}\left( \sum
 \limits_{(\boldsymbol{s},1)=n+1}2^{r_1(\boldsymbol{s},1)\theta}\right)^{\frac{1}{\theta}}\ll
 n^{-\frac{d-1}{\theta}} \left( \sum
 \limits_{(\boldsymbol{s},1)=n+1}1\right)^{\frac{1}{\theta}}\ll 1.
$$

Для $f_2$ будемо мати
$$
 \|f_2\|_{S^{\boldsymbol{r}}_{1,\infty}}\asymp \sup
 \limits_{(\boldsymbol{s},1)=n+1}2^{(\boldsymbol{s},\boldsymbol{r})}
\|A_{\boldsymbol{s}}^*(f_2,\cdot)\|_1\asymp
$$
$$
 \asymp 2^{-nr_1} \sup \limits_{(\boldsymbol{s},1)=n+1}2^{(\boldsymbol{s},\boldsymbol{r})}
\|A_{\boldsymbol{s}}^*(\cdot)\|_1\asymp 2^{-nr_1} \sup \limits_{(\boldsymbol{s},1)=n+1}2^{(\boldsymbol{s},\boldsymbol{r})}
 \ll 1.
$$

Далі, враховуючи вибір функцій $f_1$ i $f_2$, маємо
$S_{Q_n^{\gamma}}(f_1,\boldsymbol{x})=0$ i ${S_{Q_n^{\gamma}}(f_2,\boldsymbol{x})=0}$. Таким чином,
беручи до уваги оцінку $(\ref{As_Sum_infty})$, будемо мати
$$
\|f_1(\cdot)-S_{Q_n^{\gamma}}(f,\cdot)\|_{\infty}=\|f_1(\cdot)\|_{\infty}\asymp
2^{-n(r_1-1)}n^{(d-1)\left(1-\frac{1}{\theta}\right)},
$$
$$
\|f_2(\cdot)-S_{Q_n^{\gamma}}(f,\cdot)\|_{\infty}=\|f_2(\cdot)\|_{\infty}\asymp
2^{-n(r_1-1)}n^{d-1}.
$$

Оцінки знизу встановлено.

Теорему~1 доведено.

\vskip 1mm

\bf Теорема~2. \it Нехай $1<q<\infty$ i $r_1>1-\frac{1}{q}$. Тоді для $1 \leqslant
\theta \leqslant \infty$ мають місце порядкові співвідношення
\begin{equation}\label{En_teor_1q}
E_{Q_n^{\boldsymbol{\gamma}}}(S_{1,\theta}^{\boldsymbol{r}}B)_q\asymp
\mathcal{E}_{Q_n^{\boldsymbol{\gamma}}}(S_{1,\theta}^{\boldsymbol{r}}B)_q\asymp
2^{-n\left(r_1-1+\frac{1}{q}\right)}n^{(\nu-1)\left(\frac
{1}{q}-\frac {1}{\theta}\right)_+},
\end{equation}
де $a_+=\max\{a;\ 0\}$.

\vskip 2.5 mm

{\textbf{\textit{Доведення.}}} \ \rm Отримаємо спочатку оцінку зверху. Оскільки $f\in S_{1,\theta}^{\boldsymbol{r}}B$ з деяким $\boldsymbol{r}$,  ${r_1>1-\frac{1}{q}}$, то,  згідно з теоремою~A, $f\in L_q(\mathbb{R}^d)$. Тоді для ${1<q_0<q}$, скориставшись лемою~Б, а далі застосувавши нерівність різних метрик, отримуємо
 $$
 \mathcal{E}_{Q_n^{\boldsymbol{\gamma}}}(f)_q= \|f(\cdot)-S_{Q_n^{\boldsymbol{\gamma}}}(f,\cdot)\|_q = \bigg \| \sum \limits_{(\boldsymbol{s},\boldsymbol{\gamma})>n} \delta_{\boldsymbol{s}}^*(f,\cdot) \bigg \|_q \ll
 $$
 $$
 \ll\left( \sum \limits_{(\boldsymbol{s},\boldsymbol{\gamma})>n}\|\delta_{\boldsymbol{s}}^*
 (f,\cdot)\|_{q_0}^q 2^{\|s\|_1\left(\frac{1}{q_0}-\frac{1}{q}\right)q} \right)^{\frac
{1}{q}}\asymp \left( \sum \limits_{(\boldsymbol{s},\boldsymbol{\gamma})>n}\|A_{\boldsymbol{s}}^*
 (f,\cdot)\|_{q_0}^q 2^{\|s\|_1\left(\frac{1}{q_0}-\frac{1}{q}\right)q} \right)^{\frac
{1}{q}}\ll
 $$
$$
\ll \left( \sum \limits_{(\boldsymbol{s},\boldsymbol{\gamma})>n}\|A_{\boldsymbol{s}}^*
 (f,\cdot)\|_{1}^q 2^{\|s\|_1\left(1-\frac{1}{q_0}\right)q} 2^{\|s\|_1\left(\frac{1}{q_0}-\frac{1}{q}\right)q} \right)^{\frac
{1}{q}}=
$$
$$
=\left( \sum \limits_{(\boldsymbol{s},\boldsymbol{\gamma})>n}\|A_{\boldsymbol{s}}^*
 (f,\cdot)\|_{1}^q 2^{\|s\|_1\left(1-\frac{1}{q}\right)q}  \right)^{\frac
{1}{q}}=:J_2.
$$
Для того, щоб продовжити  оцінку $J_2$ розглянемо декілька випадків.

Нехай $1<q < \theta < \infty$. Тоді, застосувавши до $J_2$ нерівність
Гельдера з показником $\frac{\theta}{q}$ та врахувавши, що $r_1>1-\frac1q$, одержимо
$$
 J_2 =
 \left( \sum \limits_{(\boldsymbol{s},\boldsymbol{\gamma})>n}\|A_{\boldsymbol{s}}^*
 (f,\cdot)\|^q_1 \, 2^{(\boldsymbol{s},\boldsymbol{r})q} \, 2^{-(\boldsymbol{s},\boldsymbol{r})q} \; 2^{\|s\|_1\left(1-\frac{1}{q}\right)q}
 \right)^{\frac {1}{q}}\ll
 $$
$$
\ll \left( \sum \limits_{(\boldsymbol{s},\boldsymbol{\gamma})>n}\|A_{\boldsymbol{s}}^*
 (f,\cdot)\|^{\theta}_1 2^{(\boldsymbol{s},\boldsymbol{r})\theta}\right)^{\frac {1}{\theta}}
 \left( \sum \limits_{(\boldsymbol{s},\boldsymbol{\gamma})>n}
 \left(2^{-(\boldsymbol{s},\boldsymbol{r})q} \; 2^{\|s\|_1\left(1-\frac{1}{q}\right)q}\right)^{\frac{\theta}{\theta-q}}\right)^{\frac {1}{q}-\frac
 {1}{\theta}} \leqslant  $$
 $$
 \leqslant \|f\|_{S_{1,\theta}^{\boldsymbol{r}}B} \left( \sum \limits_{(\boldsymbol{s},\boldsymbol{\gamma})>n }
 \left(2^{-\left((\boldsymbol{s},\boldsymbol{r})-\left(1-\frac{1}{q}\right)\|s\|_1\right)}
 \right)^{\frac{q\theta}{\theta-q}}\right)^{\frac {1}{q}-\frac
 {1}{\theta}}\leqslant
 $$
$$
\leqslant \left( \sum \limits_{(\boldsymbol{s},\boldsymbol{\gamma})>n }
 \left(2^{-\left(\boldsymbol{s},\boldsymbol{r}-\left(1-\frac{1}{q}\right)
\right)} \right)^{\frac{q\theta}{\theta-q}}\right)^{\frac {1}{q}-\frac
 {1}{\theta}}=
 \left( \sum \limits_{(\boldsymbol{s},\boldsymbol{\gamma})>n }
 2^{-(\boldsymbol{s},\bar{\boldsymbol{\gamma}})\left(r_1-1+\frac{1}{q} \right)\frac{q\theta}{\theta-q}} \right)^{\frac {1}{q}-\frac
 {1}{\theta}},
$$
де $\bar{\boldsymbol{\gamma}}$ вектор з координатами $\bar{\gamma_j}=\frac{r_j-1+\frac{1}{q}}{r_j-1+\frac{1}{q}}$, $j=\overline{1,d}$. Легко переконатися, що $\bar{\gamma_j}=\gamma_j$, $j=\overline{1,\nu}$ і  $\bar{\gamma_j}\geqslant\gamma_j$, $j=\overline{\nu+1,d}$. Застосувавши до останньої суми лему~Г, отримуємо
$$
J_1\ll  2^{-n\left(r_1-1+\frac{1}{q}\right)}
 n^{(\nu-1)\left(\frac {1}{q}-\frac{1}{\theta}\right)}.
 $$

 У випадку $1\leqslant \theta \leqslant q<\infty$, $q\neq1$, скориставшись нерівністю
 $$
 \left( \sum \limits_k|a_k|^{v_2}\right)^{\frac{1}{v_2}}\leqslant
 \left( \sum \limits_k|a_k|^{v_1}\right)^{\frac{1}{v_1}}, \
 0<v_1\leqslant v_2<\infty \ ,
 $$
(див.,~\cite[~с.~43]{Xardi}), застосувавши нерівність Гельдера та беручи до уваги, що $r_1 >1-\frac{1}{q}$,
 оцінку $J_2$ можемо продовжити таким чином
 $$
 J_2 \leqslant \left( \sum \limits_{(\boldsymbol{s},\boldsymbol{\gamma})>n }\|A_{\boldsymbol{s}}^*
 (f,\cdot)\|^{\theta}_1
 2^{\|s\|_1\left(1-\frac{1}{q}\right)\theta}
 \right)^{\frac {1}{\theta}} =
 $$
 $$
 = \left( \sum \limits_{(\boldsymbol{s},\boldsymbol{\gamma})>n }\|A_{\boldsymbol{s}}^*
 (f,\cdot)\|^{\theta}_1  2^{(\boldsymbol{s},\boldsymbol{r})\theta}
 2^{-(\boldsymbol{s},\boldsymbol{\bar{\gamma}})\left(r_1-1+\frac{1}{q}\right)\theta}
 \right)^{\frac {1}{\theta}} \ll
 $$
 $$
 \ll \left( \sum \limits_{(\boldsymbol{s},\boldsymbol{\gamma})>n } 2^{(\boldsymbol{s},\boldsymbol{r})\theta} \|A_{\boldsymbol{s}}^*
 (f,\cdot)\|^{\theta}_1
 \right)^{\frac {1}{\theta}}  \sup \limits_{(\boldsymbol{s},\boldsymbol{\gamma})>n }2^{-(\boldsymbol{s},\boldsymbol{\bar{\gamma}})\left(r_1-1+\frac{1}{q}\right)}\leqslant
 $$
 $$
 \leqslant  \|f\|_{S^{\boldsymbol{r}}_{1,\theta}B} 2^{-n\left(r_1-1+\frac{1}{q}\right)}\leqslant 2^{-n\left(r_1-1+\frac{1}{q}\right)},
 $$
де, як і в попередньому випадку,  вектор $\bar{\boldsymbol{\gamma}}$ визначається аналогічно й $\bar{\boldsymbol{\gamma}}\geqslant\boldsymbol{\gamma}$.

 Нехай тепер $\theta = \infty$. Тоді для $f\in S^{\boldsymbol{r}}_{1,\infty}B$
 згідно з (\ref{Norm_dek_Sr1_inf}) маємо
 $$
 \|A_{\boldsymbol{s}}^*(f,\cdot)\|_1 \ll 2^{-(\boldsymbol{s},\boldsymbol{r})},
 $$
 а тому, застосувавши лему~Д, отримуємо
 $$
 J_1\ll \left( \sum \limits_{(\boldsymbol{s},\boldsymbol{\gamma})>n} 2^{-(\boldsymbol{s},\boldsymbol{r})q} 2^{\|s\|_1\left(1-\frac{1}{q}\right)q} \right)^{\frac
{1}{q}}=  \left( \sum \limits_{(\boldsymbol{s},\boldsymbol{\gamma})>n} 2^{-(\boldsymbol{s},\boldsymbol{\bar{\gamma}})\left(r_1-1+\frac{1}{q}\right)q} \right)^{\frac
{1}{q}}  \asymp
$$
$$
\asymp 2^{-n\left(r_1-1+\frac{1}{q}\right)} n^{\frac
 {\nu-1}{q}}.
 $$

Оцінки зверху в теоремі встановлено.

Перейдемо до встановлення оцінок знизу. Для цього при певних значеннях параметрів $q$ і $\theta$ достатньо вказати функції $f\in S^{\boldsymbol{r}}_{1,\theta}B$, для яких оцінки
знизу величин $\mathcal{E}_{Q_n^{\boldsymbol{\gamma}}}(f)_q$ співпадають за порядком з
оцінками знизу величин
$\mathcal{E}_{Q_n^{\boldsymbol{\gamma}}}(S^{\boldsymbol{r}}_{1,\theta}B)_q$ в (\ref{En_teor_1q}). Зауважимо, що достатньо розглянути випадок $\nu=d$, тобто будемо вважати $\gamma_j=1$, $j=\overline{1,d}$.

Нехай $1\leqslant \theta \leqslant q$, $q\neq 1$. Розглянемо функцію
$$
f_3(\boldsymbol{x})=2^{-r_1n}A_{\boldsymbol{\tilde{s}}}^*(\boldsymbol{x}),
$$
де $\|\boldsymbol{\tilde{s}}\|_1=n+1$.

Покажемо, що $f_3\in S_{1,\theta}^{\boldsymbol{r}}B$. Маємо
$$
\|f_3(\cdot)\|_{S_{1,\theta}^{\boldsymbol{r}}B}\asymp \left(\sum\limits_{s\geqslant 0} 2^{(\boldsymbol{s},\boldsymbol{r})\theta} \|A_{\boldsymbol{s}}^*(f_3,\cdot)\|_1^{\theta} \right)^{\frac{1}{\theta}}\asymp
$$
$$
\asymp 2^{-r_1n} \left(2^{(\boldsymbol{\tilde{s}},\boldsymbol{r})\theta} \|A_{\boldsymbol{\tilde{s}}}^*(\cdot)\|_1^{\theta} \right)^{\frac{1}{\theta}}\ll 2^{-r_1n}  2^{-r_1n} =1.
$$

Оскільки для функцій $f_3(\boldsymbol{x})$  виконується співвідношення
$S_{Q_n^{\gamma}}(f_3,\boldsymbol{x})=0$, то згідно з лемою~2, будемо мати
$$
\mathcal{E}_{Q_n^{\boldsymbol{\gamma}}}(S^{\boldsymbol{r}}_{1,\theta}B)_q\gg \|f_3(\cdot)-S_{Q_n^{\gamma}}(f_3,\cdot)\|_{q}=\|f_3(\cdot)\|_q\asymp
2^{-nr_1}\|A_{\boldsymbol{\tilde{s}}}^*(\cdot)\|_q\asymp
$$
$$
\asymp 2^{-nr_1} 2^{n\left(1-\frac{1}{q}\right)} = 2^{-n\left(r_1-1+\frac{1}{q}\right)}.
$$

У випадку $1< q < \theta <\infty$ розглянемо функцію
$$
f_4(\boldsymbol{x})=2^{-nr_1} n^{-\frac{d-1}{\theta}}\sum\limits_{\|\boldsymbol{s}\|_1=n+1}A_{\boldsymbol{s}}^*(\boldsymbol{x}).
$$

Покажемо, що $f_4\in S_{1,\theta}^{\boldsymbol{r}}B$. Маємо
$$
\|f_4(\cdot)\|_{S_{1,\theta}^{\boldsymbol{r}}B}\asymp \left(\sum\limits_{s\geqslant 0} 2^{(\boldsymbol{s},\boldsymbol{r})\theta} \|A_{\boldsymbol{s}}^*(f_4,\cdot)\|_1^{\theta} \right)^{\frac{1}{\theta}}=
$$
$$
=2^{-nr_1} n^{-\frac{d-1}{\theta}} \left(\sum\limits_{\|\boldsymbol{s}\|_1=n+1} 2^{(\boldsymbol{s},\boldsymbol{r})\theta} \|A_{\boldsymbol{s}}^*(\cdot)\|_1^{\theta}\right)^{\frac{1}{\theta}}\ll 2^{-nr_1} n^{-\frac{d-1}{\theta}} \left(\sum\limits_{\|\boldsymbol{s}\|_1=n+1} 2^{(\boldsymbol{s},\boldsymbol{r})\theta} \right)^{\frac{1}{\theta}} \ll 1.
$$

Врахувавши, що за рахунок вибору функцій $f_4$  виконується співвідношення
$S_{Q_n^{\gamma}}(f_4,\boldsymbol{x})=0$,  будемо мати
$$
\mathcal{E}_{Q_n^{\boldsymbol{\gamma}}}(S^{\boldsymbol{r}}_{1,\theta}B)_q\gg \|f_4(\cdot)-S_{Q_n^{\gamma}}(f_4,\cdot)\|_{q}=\|f_4(\cdot)\|_q.
$$

Оскільки, як було показано вище $f_4\in S^{\boldsymbol{r}}_{1,\theta}B$, а за умовами теореми ${r_1>1-\frac{1}{q}}$, то, згідно з теоремою~A, $f_4\in L_q(\mathbb{R}^d)$. Для $\boldsymbol{s}\in \mathbb{Z}^d_+$ покладемо
 $$
 \Delta(\boldsymbol{s})=\Big\{\boldsymbol{x}: 2^{-s_j-1}\leqslant x_j<2^{-s_j},\ j=\overline {1,d}
 \Big\},
 $$
 $\Delta(\boldsymbol{s}) \cap \Delta(\boldsymbol{s}')=\varnothing$, якщо $\boldsymbol{s}\neq
\boldsymbol{s}'$, тоді за теоремою~В (Літлвуда--Пелі) будемо мати
 $$
    \|f_{4}(\cdot)\|_q \gg \Bigg \| \left(
  \sum \limits_{\|\boldsymbol{s}\|_1=n+1}|\delta^*_{\boldsymbol{s}}(f_{4},\cdot)|^2
  \right)^{\frac{1}{2}} \Bigg \|_q \geqslant
 $$
 $$
  \geqslant \left(\sum \limits_{\|\boldsymbol{s}\|_1=n+1}
  \int \limits_{\Delta(\boldsymbol{s})}|\delta^*_{\boldsymbol{s}}(f_{4},\boldsymbol{x})|^q d\boldsymbol{x}\right)^{\frac{1}{q}}.
 $$

Скориставшись оцінкою (\ref{As_int}) та лемою~Д, останню оцінку можемо продовжити таким чином
$$
    \|f_{4}(\cdot)\|_q \gg  2^{-nr_1} n^{-\frac{d-1}{\theta}} \left(\sum \limits_{\|\boldsymbol{s}\|_1=n+1}
  \int \limits_{\Delta(\boldsymbol{s})}|A^*_{\boldsymbol{s}}(\boldsymbol{x})|^q d\boldsymbol{x}\right)^{\frac{1}{q}}\gg
 $$
 $$
\gg 2^{-nr_1} n^{-\frac{d-1}{\theta}}  \left(\sum \limits_{\|\boldsymbol{s}\|_1=n+1}
  2^{\|\boldsymbol{s}\|_1(q-1)}\right)^{\frac{1}{q}}\asymp
 $$
$$
\asymp 2^{-nr_1} n^{-\frac{d-1}{\theta}}  2^{n\frac{q-1}{q}} n^{\frac{d-1}{q}}= 2^{-n\left(r_1-\left(1-\frac{1}{q}\right)\right)} n^{(d-1)\left(\frac{1}{q}-\frac{1}{\theta}\right)}
$$

Насамкінець розглянемо випадок $\theta=\infty$  і, відповідно, функцію
$$
f_5(\boldsymbol{x})=2^{-nr_1} \sum\limits_{\|\boldsymbol{s}\|_1=n+1}A_{\boldsymbol{s}}^*(\boldsymbol{x}).
$$

Покажемо, що $f_5\in S_{1,\infty}^{\boldsymbol{r}}B$. Маємо
$$
\|f_5(\cdot)\|_{S_{1,\infty}^{\boldsymbol{r}}B}\asymp \sup\limits_{s\geqslant 0} 2^{(\boldsymbol{s},\boldsymbol{r})} \|A_{\boldsymbol{s}}^*(f_5,\cdot)\|_1 =2^{-nr_1} \sup\limits_{\|\boldsymbol{s}\|_1=n+1} 2^{(\boldsymbol{s},\boldsymbol{r})} \|A_{\boldsymbol{s}}^*(\cdot)\|_1 \ll 1.
$$

Врахувавши, що для функцій $f_5$  виконується співвідношення
$S_{Q_n^{\gamma}}(f_5,\boldsymbol{x})=0$, аналогічно до того, як і в попередньому випадку, отримаємо
$$
\mathcal{E}_{Q_n^{\boldsymbol{\gamma}}}(S^{\boldsymbol{r}}_{1,\theta}B)_q\gg \|f_5(\cdot)-S_{Q_n^{\gamma}}(f_5,\cdot)\|_{q}=\|f_5(\cdot)\|_q\gg
$$
$$
 \gg \Bigg \| \left(
  \sum \limits_{\|\boldsymbol{s}\|_1=n+1}|\delta^*_{\boldsymbol{s}}(f_{5},\cdot)|^2
  \right)^{\frac{1}{2}} \Bigg \|_q\geqslant \left(\sum \limits_{\|\boldsymbol{s}\|_1=n+1}
  \int \limits_{\Delta(\boldsymbol{s})}|\delta^*_{\boldsymbol{s}}(f_{5},\boldsymbol{x})|^q d\boldsymbol{x}\right)^{\frac{1}{q}}\asymp
 $$
$$
 \asymp  2^{-nr_1} \left(\sum \limits_{\|\boldsymbol{s}\|_1=n+1}
  \int \limits_{\Delta(\boldsymbol{s})}|A^*_{\boldsymbol{s}}(\boldsymbol{x})|^q d\boldsymbol{x}\right)^{\frac{1}{q}}\gg 2^{-nr_1} \left(\sum \limits_{\|\boldsymbol{s}\|_1=n+1}
  2^{\|\boldsymbol{s}\|_1(q-1)}\right)^{\frac{1}{q}}\asymp
 $$
$$
\asymp 2^{-nr_1} 2^{n\frac{q-1}{q}} n^{\frac{d-1}{q}}= 2^{-n\left(r_1-1+\frac{1}{q}\right)} n^{\frac{d-1}{q}}.
$$

Оцінки знизу встановлено. Теорему~2 доведено.

\vskip 2 mm

На завершення роботи зробимо деякі коментарі, щодо одержаних результатів.

Точні за порядком оцінки величини (\ref{EQN}) для класів Нікольського  $S^{\boldsymbol{r}}_{1}H(\mathbb{R}^d)$ в метриці простору $L_q$ при $1\leqslant q<\infty$ (теорема~2) встановлено Wang Heping та Sun Youngsheng~\cite{WangHeping_SunYongsheng_1999_AppT}. Зауважимо, що методи, які використовувалися для встановлення оцінок у теоремі~2 при $\theta=\infty$ дещо відрізняються від методів, які застосовували Wang Heping та Sun Youngsheng.  Результат теореми~1 є новим і для класів Нікольського, тобто  у випадку $\theta=\infty$.

Знаходженню точних за порядком оцінок величин (\ref{EQN}) та (\ref{EQN1}) у східчастому гіперболічному хресті для ряду інших значень параметрів $p$, $\theta$ та $q$ присвячені роботи~\cite{WangHeping_SunYongsheng_1995}, \cite{Yanchenko_YMG_2013}, \cite{Yanchenko_Zb_2013}. Дослідження класів $S^{\boldsymbol{r}}_{p,\theta}B(\mathbb{R}^d)$ з точки зору знаходження оцінок інших апроксимативних характеристик проводилися, зокрема, у роботах Wang Heping~\cite{WangHeping_1997_Q}, \cite{WangHeping_2004}, Hans-Jurgen Schmeisser, Winfried Sickel~\cite{Winfried Sickel_2004} та роботі~\cite{Yanchenko_YMG_2010_8}.

Зауважимо, що більш інтенсивно досліджуються класи Нікольського--Бєсова періодичних функцій як однієї так і багатьох змінних. Так, порядкові оцінки наближення функцій з даних класів за допомогою тригонометричних поліномів з номерами гармонік із східчастого гіперболічного хреста встановлювалися в роботах В.\,М.~Темлякова~\cite{Temlyakov_1986m}, \cite{Temlyakov_1988}, А.\,С.~Романюка~\cite{Romanyuk_91}, \cite{Romanyuk_2004msb}, \cite{Romanyuk_2008mzam}, Є\,М.~Галєєва~\cite{Galeev_90_Zam}. Більш детально з дослідженнями  класів Нікольського--Бєсова періодичних функцій, з точки зору знаходження порядкових оцінок різних апроксимативних характеристик, можна ознайомитися в монографії А.\,С.~Романюка~\cite{Romanyuk_2012m}.

На даний час також інтенсивно досліджуються й узагальнення класів Нікольського--Бєсова з домінуючою мішаною похідною. У цьому напрямі відзначимо роботи М.\,М.~Пустовойтова~\cite{Pustovoitov_94}, \cite{Pustovoitov_99}, Sun Yongsheng, Wang Heping~\cite{SunYongsheng_WangHeping_1997}, Н.\,В.~Дерев'янко~\cite{Derevyanko_2014_UMG_5},  С.\,А.~Стасюка~\cite{Stasjuk_04}, С.\,А.~Стасюка і С.\,Я.~Янченка~\cite{Stasuk_Yanchenko_Anal_math}, С.\,Я.~Янченка~\cite{Yanchenko_YMG_2010_1}.

\vskip 3 mm

\textbf{Contact information:}
Department of the Theory of Functions, Institute of Mathematics of the National
Academy of Sciences of Ukraine, 3, Tereshenkivska st., 01601, Kyiv, Ukraine.

\vskip 3 mm

E-mail: \href{mailto:Yan.Sergiy@gmail.com}{Yan.Sergiy@gmail.com}


\begin{thebibliography}{10}


\bibitem{Nikolsky_63} {\it Никольский~С.\,М.\/} Функции с доминирующей смешанной производной,
 удовлетворяющей кратному условию  Гельдера~// Сиб. мат. журн.~--- 1963.~--- \textbf{4}, №6.~---
C.~1342\,--\,1364.

\bibitem{Amanov_1965} {\it Аманов~Т.\,И.  \/}
Теоремы представления и вложения для функциональных пространств
$S^{(r)}_{p,\theta}B(\mathbb{R}_n)$ и $S^{(r)_*}_{p,\theta}B$,
($0\leqslant x_j\leqslant 2\pi$; $j=1,\ldots,n$)~// Тр. Мат. ин-та
АН СССР.~--- 1965.~--- \textbf{77}.~--- С.~5\,--\,34.


\bibitem{Lizorkin_Nikolsky_1989} {\it Лизоркин~П.\,И., Никольский~С.\,М.\/}
Пространства функций смешанной гладкости с декомпозиционной точки
зрения~// Тр. Мат. ин-та АН СССР.~--- 1989.~--- \textbf{187}.~---
C.~143\,--\,161.



\bibitem{Lizorkin_69} {\it Лизоркин~П.\,И.\/}
Обобщенное лиувиллевское дифференцирование и метод мультипликаторов
в теории вложений классов дифференцируемых функций~// Тр. Мат. ин-та
АН СССР.~--- 1969.~--- \textbf{105}.~--- C.~89\,--\,167.


\bibitem{WangHeping_SunYongsheng_1999_AppT} {\it Wang Heping, Sun Yongsheng. \/}
Approximation of functions in  $\widetilde{S_1^{r}L}$, $S_1^{r}H$ by entire functions~// Approx. Theory and its Appl.~--- 1999.~---
\textbf{11}, №4.~--- P.~88~--~93.


\bibitem{WangHeping_1997_Q} {\it Wang Heping.\/}
Quadrature formulas for classes of functions with bounded mixed derivative or difference~// Science in China (Series A).~--- 1997.~---
\textbf{40}, №5.~--- P.~449\,--\,458.



\bibitem{Nikolsky_1969_book} {\it Никольский~С.\,М.\/} Приближение
функций многих переменных и теоремы вложения.~--- М.: Наука,
1969.~--- 480~c.

\bibitem{WangHeping_SunYongsheng_1995} {\it Wang Heping, Sun Yongsheng. \/}
Approximation of multivariate functions with certain mixed
smoothness by entire functions~// Northeast. Math. J.~--- 1995.~---
\textbf{11}, №\,4.~--- P.~454\,--\,466.


\bibitem{Temlyakov_1986m} {\it Темляков В.\,Н.\/} Приближение
функций с ограниченной смешанной производной~// Тр. Мат. ин-та АН
СССР.~--- 1986.~--- {\bf 178}.~--- C.~1\,--\,112.

\bibitem{Xardi} {\it Харди~Г.\,Г., Литтльвуд~Дж.\,Е., Полиa~Г.\/} Неравенства.~--- М.: Изд-во иностр. лит., 1948.~--- 456~c.


\bibitem{Yanchenko_YMG_2013} {\it Янченко~С.~Я.\/} Наближення функцій з класів
$S^r_{p,\theta}B$ у рівномірній метриці~// Укр. мат. журн.~---
2013.~--- \textbf{65}, №\,5.~--- С.~698~--~705.

\bibitem{Yanchenko_Zb_2013} { \it Янченко~С.\,Я.\/} Оцінки апроксимативних
характеристик класів функцій $S^r_{p,\theta}B(\mathbb{R}^d)$ у
рівномірній метриці~// Теорія наближення функцій та суміжні питання:
Зб. праць Ін-ту математики НАН України.~--- 2013.~--- \textbf{10},
№\,1.~--- С.~328\,--\,340.


\bibitem{WangHeping_2004} {\it Wang Heping. \/}
Representation and approximation of multivariate function with
bounded mixed smoothness by hyperbolic wavelets~// J.
Math. Anal. Appl.~--- 2004.~--- \textbf{291}.~--- P.~698\,--\,715.


\bibitem{Winfried Sickel_2004} {\it Hans-Jurgen Schmeisser, Winfried Sickel. \/}
Spaces of functions of mixed smoothness and approximation from
hyperbolic crosses~// J.
of Approx. Theory.~--- 2004.~--- V.~128.~--- P.~115\,--\,150.


\bibitem{Yanchenko_YMG_2010_8} {\it Янченко~С.\,Я.\/} Наближення класів
$S^{r}_{p,\theta}B(\mathbb{R}^d)$ функцій багатьох змінних цілими
функціями спеціального вигляду~// Укр. мат. журн.~--- 2010.~---
\textbf{62}, №8.~--- С.~1124\,--\,1138.


\bibitem{Temlyakov_1988} {\it Темляков~В.~Н.\/} Оценки асимптотических характеристик
классов функций с ограниченной смешанной производной или
разностью~// Тр. Мат. ин-та АН СССР.~--- 1988.~--- \textbf{189}.~---
C.~138\,--\,168.


\bibitem{Romanyuk_91} {\it Романюк А.\,С. \/}
Приближение классов Бесова периодических функций многих переменных в
пространстве $L_q$~// Укр. мат. журн.~--- 1991.~--- {\bf 43},
№10.~--- С.~1398\,--\,1408.

\bibitem{Romanyuk_2004msb} {\it Романюк~А.~С.\/} Приближение классов $B^r_{p,\theta}$
периодических функций многих переменных линейными методами и
наилучшие приближения~// Мат. сб.~--- 2004.~--- \textbf{195},
№2.~--- С.~91\,--\,116.


\bibitem{Romanyuk_2008mzam} {\it Романюк~А.~С.\/} Наилучшие тригонометрические приближения
классов периодических функций многих переменных в равномерной
метрике~// Мат. заметки.~--- 2007.~--- \textbf{82}, 2.~---
С.~247\,--\,261.

\bibitem{Galeev_90_Zam} {\it Галеев~Э.~М.\/} Приближение классов периодических функций
нескольких переменных ядерными операторами~// Мат.
заметки.~--- 1990.~--- \textbf{47}, №\,3.~--- С.~32\,--\,41.

\bibitem{Romanyuk_2012m} {\it Романюк~А.~С.\/} Аппроксимативные характеристики классов периодических функций многих переменних~// Праці Інституту математики НАН України.~--- 2012.~--- \textbf{93}.~--- 352~с.


\bibitem{Pustovoitov_94} {Пустовойтов~Н.~Н. \/} Представление и
приближение периодических функций многих переменных с заданным
смешанным модулем непрерывности~// Anal.
Math.~--- 1994.~--- Vol.~20.~--- P.~35\,--\,48.

\bibitem{Pustovoitov_99} {\it Пустовойтов~Н.\,Н.} Приближение многомерных функций с заданной
мажорантой смешанных модулей непрерывности~// Мат. заметки.~--- 1999.~--- {\bf65}, №1.~--- C.~107\,--\,117.


\bibitem{SunYongsheng_WangHeping_1997} {\it Sun Yongsheng, Wang
Heping. \/}
Representation and approximation of multivariate periodic functions
with bounded mixed moduli of smoothness~// Тр. Мат. ин-та РАН.~--- 1997.~--- T.~219.~---
С.~356\,--\,377.

\bibitem{Derevyanko_2014_UMG_5} {\it Дерев'янко~Н.\,В.\/} Наближення класів $H^{\Omega}_{p}$  періодичних функцій багатьох змінних у просторі $L_p$~// Укр. мат. журн.~---
2014.~--- \textbf{66}, №\,5.~--- С.~634\,--\,644.

\bibitem{Stasjuk_04} {\it Стасюк~С.~А.} Найкращі наближення, колмогоровські
та тригонометричні поперечники класів ${B^{\Omega}_{p,\theta}}$
періодичних функцій багатьох змінних~// Укр. мат.
журн.~--- 2004.~--- Т.~56, №~11.~--- C.~1557~--~1568.



\bibitem{Stasuk_Yanchenko_Anal_math} {\it Stasyuk~S.~А., Yachenko~S.~Ya.\/} Approximation of functions from Nikolskii--Besov type classes of generalized mixed smoothness~// Anal. Math.~--- 2015.~--- \textbf{41}.~---
P.~311\,--\,334.

\bibitem{Yanchenko_YMG_2010_1} { \it Янченко~С.\,Я.\/} Наближення класів
$B^{\Omega}_{p,\theta}$ функцій багатьох змінних цілими функціями у просторі
$L_q(\mathbb{R}^d)$~// Укр. мат. журн.~--- 2010.~--- \textbf{62},
№\,1.~--- С.~123\,--\,135.

\end{thebibliography}
\end{document}